\documentclass[12pt, a4paper]{article}
\usepackage[pagewise]{lineno}
\usepackage{graphicx}
\usepackage{amssymb}
\usepackage{latexsym, bm}
\usepackage{multicol}
\usepackage{indentfirst}
\usepackage{amssymb,amsfonts}
\usepackage{amsmath}
\usepackage{cases}
\usepackage[noadjust]{cite}
\usepackage[colorlinks,citecolor=red]{hyperref}
\newtheorem{theorem}{Theorem}
\newtheorem{lemma}{Lemma}

\newtheorem{definition}{Definition}

\newtheorem{remark}{Remark}

\usepackage{amssymb}
\numberwithin{equation}{section}
\numberwithin{theorem}{section}
\numberwithin{remark}{section}
\numberwithin{definition}{section}
\numberwithin{lemma}{section}
\numberwithin{corollary}{section}
\numberwithin{proposition}{section}
\title{The localized characterization for the singularity formation in the Navier-Stokes equations}
\author{Wenke Tan\footnote{tanwenkeybfq@163.com}\\
{\small Key Laboratory of Computing and Stochastic Mathematics (Ministry of Education),}\\
{\small School of Mathematics and Statistics, Hunan Normal University,}\\
{\small Changsha, Hunan 410081, China}\\
}

\date{}

\begin{document}
\maketitle
{\bf Abstract:}
This paper is concerned with the localized behaviors of the solution $u$ to the Navier-Stokes equations near the potential singular points. We establish the concentration rate for the $L^{p,\infty}$ norm of $u$ with $3\leq p\leq\infty$. Namely, we show that if $z_0=(t_0,x_0)$ is a singular point, then for any $r>0$, it holds
\begin{align}
\limsup_{t\to t_0^-}||u(t,x)-u(t)_{x_0,r}||_{L^{3,\infty}(B_r(x_0))}>\delta^*,\notag
\end{align} and
 \begin{align}
 \limsup_{t\to t_0^-}(t_0-t)^{\frac{1}{\mu}}r^{\frac{2}{\nu}-\frac{3}{p}}||u(t)||_{L^{p,\infty}(B_r(x_0))}>\delta^*\notag\\
 for~3<p\leq\infty, ~\frac{1}{\mu}+\frac{1}{\nu}=\frac{1}{2}~and~2\leq\nu\leq\frac{2}{3}p,\notag
 \end{align}where $\delta^*$ is a positive constant independent of $p$ and  $\nu$. Our main tools are some $\varepsilon$-regularity criteria in $L^{p,\infty}$ spaces and an embedding theorem from $L^{p,\infty}$ space into a Morrey type space. These are of independent interests.

\medskip

\section{Introduction}
We consider the concentration phenomenon near the potential singularity for the three dimensional  incompressible Navier-Stokes equations
\begin{equation}\label{NS}
 \left\{\begin{array}{ll}
\partial_tu-\Delta u+u\cdot\nabla u+\nabla P=0,\\
\nabla\cdot u=0,\\
u(0,x)=u_0(x)
\end{array}\right.
\end{equation}
where the unknowns $u$, $P$ denote the velocity vector field, pressure respectively.

It is well-known that if $u_0$ is smooth enough, then problems \eqref{NS} have a unique regular solution on $[0,T)$ for some $T>0$; see, for example, \cite{Leray,Fujita,Kaniel,Kiselev,Serrin,Wahl} and the references therein. The global existence and regularity problem of the Navier-Stokes equations is one of the most significant open questions in the field of partial differential equations. The case of $\mathbb{R}^2$ was known to Leray \cite{Leray1} in 1933. Later, the case of 2D domains with boundary was settled by Ladyzhenskaya \cite{Lady} in 1959. In the case $n=3$, some remarkable progress has been made since the pioneering work by Leray in the 1930s.  The fundamental papers of Leray \cite{Leray} and Hopf \cite{Hopf} showed the global existence of weak solutions in the whole space and on bounded open domain with smooth boundary respectively. The weak solutions, called Leray-Hopf weak solutions, satisfy \eqref{NS} in the distributional sense and belong to $L^\infty L^2\cap L^2H^1$. Moreover, the following strong global energy inequality holds
\begin{align}\label{1.2}
||u(t)||^2_{L^2}+2\int_{t_0}^{t}\int_{\mathbb{R}^3}|\nabla u|^2dxdt\leq||u_(t_0)||^2_{L^2},
\end{align}
for all $t\in(0,\infty)$ and a.e. $t_0\in[0,t]$ including $0$. The regularity or uniqueness of Leray-Hopf weak solutions is one of the most significant open questions in the field of partial differential equations.

To understand the weak solutions of the Navier-Stokes equations in dimension $n=3$, there are various sufficient conditions to ensure the regularity of weak solutions.
Leray \cite{Leray} shown that for $3<p\leq\infty$, there exists $c_p$ such that the conditions
\begin{align}\label{1.3}
||u(t)||_{L^p}<\frac{c_p}{(T-t)^\frac{p-3}{2p}}
\end{align} imply the regularity of weak solutions on $[0,T]$. 
The well-known Ladyzhenskaya-Prodi-Serrin criteria \cite{Lady2,Prodi,Serrin1} showed that if $u\in L^q([0,T],L^p(\mathbb{R}^3))$ for $\frac{2}{q}+\frac{3}{p}\leq1,p>3$ then $u$ is regular on $[0,T]$. The endpoint case $ p=3$ is more subtle. In a breakthrough paper, Escauriaza, Serengin and Sverak \cite{Escauriaza} proved that the $L^\infty L^3$ solutions are smooth. This result was improved by Tao \cite{Tao} showed that as the solution $u$ approaches a finite blowup time $T$, the critical norm $||u(t)||_{L^3}$ must blow up at a rate $(\log\log\log\frac{1}{T-t})^{c}$ with some absolute constant $c>0$. The other endpoint case $p=\infty$ was generalized by Kozono and Taniuchi. In \cite{Kozono}, they proved that $u\in L^2([0,T];BMO(\mathbb{R}^3))$ implies the regularity of the solution $u$ to \eqref{NS}. Since the condition $||u(t)||_{L^p(\mathbb{R}^3)}\leq\frac{c_p}{|T-t|^\frac{p-3}{2p}}$ merely implies $u\in L^{q,\infty}(L^p)$ for $\frac{3}{p}+\frac{2}{q}=1,p>3$, it is natural to generalize the classical Ladyzhenskaya-Prodi-Serrin type criterion in Lorentz spaces. In \cite{Kozono1}, Kim and Kozono proved the local boundedness of a weak solution $u$ under the assumption that $||u||_{L^{r,\infty}([0,T];L^{s,\infty}(\mathbb{R}^3))}$ is sufficiently small for some $(r,s)$ with $\frac{2}{r}+\frac{3}{s}=1$ and $3\leq s<\infty$. The limiting case of the regularity criteria derived by Kim and Kozono was proved by He and Wang \cite{He} i.e. any weak solution $u$ to the Navier-Stokes equations is regular under the assumption that $||u||_{L^{2,\infty}([0,T];L^\infty(\mathbb{R}^3))}$ is sufficiently small. This results of He and Wang were improved by Wang and Zhang \cite{Wang} which showed that $||u_3||_{L^{r,\infty}([0,T];L^{s,\infty}(\mathbb{R}^3))}\leq M$ and $||u_h||_{L^{r,\infty}([0,T];L^{s,\infty}(\mathbb{R}^3))}\leq c_M$ with $\frac{2}{r}+\frac{3}{s}=1$ and $3<s\leq\infty$ imply the regularity of the suitable weak solution $u$ to Navier-Stokes equations, where $c_M$ is a small constant depending on $M$. 

Another important step towards a better understanding of the Navier-Stokes equations is the partial regularity theory. This theory was initiated by Scheffer \cite{Sch1,Sch2,Sch3} and improved by Caffarelli, Kohn and Nirenberg \cite{Caffarelli}. In \cite{Sch3}, Scheffer pioneered the partial regularity theory by introducing the definition of suitable weak solutions and proving their existence in dimension $n=3$. Moreover, he showed that the singular sets of the suitable weak solutions have finite $\frac{5}{3}$-dimensional Hausdorff measure in space-time. Caffarelli, Kohn and Nirenberg \cite{Caffarelli} made remarkable improvements in dimension $n=3$ by proving that the $1$-dimensional parabolic Hausdorff measure of singular sets of suitable weak solutions is zero.  For more results about partial regularity theory of the Navier-Stokes equations, we refer the reader to \cite{C-L,G-K-T,Ladyzhenskaya,Lin,Wolf,W-W-Z} and the references therein.

On the other hand, the idea of investigating the potential singularity of solutions goes back as far as \cite{Leray}. In \cite{Leray}, Leray showed that if a weak solution $u$ first develops singularity at time $T$ then for $3<p\leq\infty$ and $t<T$, it follows
\begin{align}\label{1.4}
||u(t)||_{L^p(\mathbb{R}^3)}\geq\frac{c_p}{(T-t)^\frac{p-3}{2p}}.
\end{align} Moreover, Leray raised the question of the existence of self-similar singularity with the form
\begin{align*}
u(x,t)=\frac{1}{\sqrt{2a(T-t)}}U(\frac{x}{\sqrt{2a(T-t)}}).
\end{align*} This question was completely solved by a negative answer due to Ne\v{c}as, Ru\v{z}i\v{c}ka and \v{S}ver\'{a}k \cite{N-R-S}, also see \cite{Tsai} for a more general case. In general, if $u$ satisfies
\begin{align*}
||u(t)||_{L^p(\mathbb{R}^3)}\leq\frac{C}{(T-t)^\frac{p-3}{2p}},\end{align*} The singularity or regularity of solution $u$ at time $T$ remains unknown.  The potential singularity satisfying
\begin{align*}
||u(t)||_{L^\infty(\mathbb{R}^3)}\leq\frac{C}{(T-t)^\frac{1}{2}}\end{align*} is called Type I singularity in time. For the axisymmetric Navier-Stokes equations, Chen-Strain-Yau-Tsai \cite{C-S-Y-T1,C-S-Y-T2} and Koch-Nadirashvili-Seregin-\v{S}ver\'{a}k \cite{K-N-S-S} proved that the solution $u$ does not develop Type I singularity respectively.
For the behavior of the critical $L^3$ norm, Escauriaza, Seregin and Sver\'{a}k \cite{Escauriaza} proved that if $(x,T)$ is a singular point then
\begin{align}\label{1.5}
\limsup_{t\to T^-}||u(t)||_{L^3(B_r(x))}=\infty~for~any~fixed~r>0.
\end{align} Later, Seregin \cite{Seregin2} improved \eqref{1.5}:
\begin{align}\label{1.6}
\lim_{t\to T^-}||u(t)||_{L^3(\mathbb{R}^3)}=\infty.
\end{align} Albritton and Barker \cite{A-B} refined \eqref{1.5} and \eqref{1.6} to show that if $\Omega$ is a bounded domain with $C^2$ boundary one has
\begin{align}\label{1.7}
\lim_{t\to T^-}||u(t)||_{L^3(B_\delta(x)\cap\Omega)}=\infty~for~any~fixed~\delta>0.
\end{align}In \cite{L-O-W}, Li, Ozawa and Wang proved that if $u$ first blows up at $T$, there exists $t_n\to T^-$ and $x_n\in\mathbb{R}^3$ such that
\begin{align}\label{1.8}
||u(t_n)||_{L^m(B_{\sqrt{C(m)(T-t)}})(x_n)}\geq\frac{C(m)}{(T-t)^\frac{m-3}{2m}}~for~3\leq m\leq\infty.
\end{align} This result was improved by Maekawa, Miura and Prange. They \cite{M-M-P} proved that for every $t\in(0,T)$ there esists $x(t)\in\mathbb{R}^3$ such that
\begin{align}\label{1.9}
||u(t)||_{L^m(B_{\sqrt{C(m)(T-t)}(x(t))})}\geq\frac{C(m)}{(T-t)^\frac{m-3}{2m}}~for~3\leq m\leq\infty.
\end{align}It is worth pointing out that in \eqref{1.8} and \eqref{1.9}, there is no information about $x_n$ and $x(t)$. It is natural to ask whether the concentration phenomenon occurs on balls $B(x,R)$ with $R=O(\sqrt{T-t})$ and with $(x,T)$ being a singular point. Recently, this question was affirmatively answered by Barker and Prange for the critical $L^3$ norm for Leray-Hopf solutions which experiences the first singular time at $T$. In \cite{B-P}, they proved that if $u$ satisfy the Type I bound:
\begin{align}\label{1.10}
\sup_{\bar{x}\in\mathbb{R}^3}\sup_{0<r<r_0}\sup_{T-r^2<t<T}(\frac{1}{r}\int_{B_r(\bar{x})}|u(y,t)|^2dy)^\frac{1}{2}\leq M\\~for~a~given~r_0\in(0,\infty]~and~M,~T\in(0,\infty)\notag
\end{align}then it holds
\begin{align}\label{1.11}
||u(\cdot,t)||_{L^3(B_R(x))}\geq\gamma_{univ},\quad R=O(\sqrt{T-t}).
\end{align} Recently, Barker and Prange \cite{B-P1} show under the assumption $||u||_{L^\infty_t L^{3,\infty}}\leq M$, the optimal blow-up rate at the potential singular point $(T^*,0)$ is
\begin{align}
||u(t,\cdot)||_{L^3(B_R(0))}\geq C(M)\log(\frac{1}{T^*-t}),~R=O((T^*-t)^{\frac{1}{2}-}).
\end{align}

In conclusion, if one characterizes the local behaviors of solutions to the Navier-Stokes equations near a potential singularity by critical norm $L^3$, the optimal blow-up rate was obtained by Barker and Prange \cite{B-P1}. But, if we consider the local characterization of singularity by $L^3$ norm, there is an unpleasant problem. On the one hand, It is well-known that if $|u(T,x)|\leq\frac{c}{|x|}$ with small enough $c$ then $(T,0)$ can not be a singular point. One the other hand, it is clear that $||u(T,\cdot)||_{L^3(B_r)}=\infty$ for any $r>0$. This means that one can not exclude such point from the singular set by using $L^3$ norm. Noticing that $||u(T,\cdot)||_{L^{3,\infty}(B_r)}=(\frac{4\pi}{3})^\frac{1}{3}c$, it is more natural to characterize the singularity formation for the Navier-Stokes equations by $L^{3,\infty}$ norm. The main purpose of this paper is to characterize the singularity formation in the Navier-Stokes equations by the critical norm $L^{3,\infty}$.
\subsection{Main result}
We first show some $\varepsilon$-regularity criteria. It is worth pointing out that our criteria are established in $L^{q,\infty}$ space and the constant $\delta$ in our $\varepsilon$-regularity criteria does not depend on the value of $p$.
\begin{theorem}\label{mainr}
Let $\frac{1}{q}+\frac{1}{p}=\frac{1}{2}$ with $2\leq p\leq\infty$. Assume $(u,P)$ be a suitable weak solution to the Navier-Stokes equations \eqref{NS} on $Q_1(z_0)$. There exists an absolute constant $\delta>0$ such that if
\begin{align}
||\sup_{\eta\leq 1}(\frac{1}{\eta}\int_{B_\eta(x_0)}|u(x,t)-u_{x_0,\eta}|^pdx)^\frac{1}{p}||_{L^{q,\infty}[t_0-1,t_0]}\leq&\delta\\
or~~~~~~~ ||\sup_{\eta\leq 1}(\frac{1}{\eta}\int_{B_\eta(x_0)}|u(x,t)|^pdx)^\frac{1}{p}||_{L^{q,\infty}[t_0-1,t_0]}\leq&\delta,
\end{align}
 then $z_0$ is a regular point.
\end{theorem}
\begin{remark}
It is worth pointing out that the quantities
\begin{align*}
&||(\sup_{\eta\leq 1}\frac{1}{\eta}\int_{B_\eta(x_0)}|u(x,t)-u_{x_0,\eta}|^pdx)^\frac{1}{p}||_{L^{q,\infty}[t_0-1,t_0]}\\ and~ &||\sup_{\eta\leq 1}(\frac{1}{\eta}\int_{B_\eta(x_0)}|u(x,t)|^pdx)^\frac{1}{p}||_{L^{q,\infty}[t_0-1,t_0]}
\end{align*} with $\frac{1}{q}+\frac{1}{p}=\frac{1}{2}$ and $2\leq p\leq\infty$ are invariant under the scaling \eqref{1.17}, we can replace $1$ by any $r>0$ in Theorem \ref{mainr}.
\end{remark}

By using Theorem \ref{mainr} and the embedding theorem established in Lemma \ref{l4}, we obtain the following theorem.
\begin{theorem}\label{mains}
Let $(u,P)$ be a suitable weak solution in $Q_1(z_0)$. Assume $z_0$ be a singular point. Then
for any given $r\in(0,1)$, it holds
\begin{align}
&\limsup_{t\to t_0^-}||u(t,x)-u(t)_{x_0,r}||_{L^{3,\infty}(B_r(x_0))}>\delta^*\\ and\notag\\
&\limsup_{t\to t_0^-}(t_0-t)^{\frac{1}{\mu}}r^{\frac{2}{\nu}-\frac{3}{p}}||u(t)||_{L^{p,\infty}(B_r(x_0))}>\delta^*\\&for~3<p\leq\infty, ~\frac{1}{\mu}+\frac{1}{\nu}=\frac{1}{2}~and~2\leq\nu\leq\frac{2}{3}p,\notag
\end{align} where $u_{x_0,r}(t)=\frac{1}{|B_r(x_0)|}\int_{B_r(x_0)}u(t,y)dy$ and $\delta^*>0$ is independent on $\mu,\nu,p$ and $r$.
\end{theorem}

 Before the proofs of main results, we first recall some definitions and notations of the suitable weak solutions to \eqref{NS}, Lorenz space, and some invariant quantities.
 Setting
\begin{align*}
B_r(x_0)=\{x\in\mathbb{R}^3:|x-x_0|<r\},B_r=B_r(0),~B=B_1,\\
Q_r(z_0)=B_r(x_0)\times(t_0-r^2,t_0),~Q_r=Q_r(0),~Q=Q_1.
\end{align*}

\begin{definition}
The function pair $(u,P)$ is called a suitable weak solution of \eqref{NS} in
$Q_1(z_0)$ if\\
1.$u\in L^\infty((t_0-1,t_0); L^2_{loc}(B_1(x_0))\cap L^2((t_0-1,t_0); H^1_{loc}(B_1(x_0))$,\\
2.There exists a distribution $P\in L^\frac{3}{2}_{loc}(Q_1(z_0))$ such that $(u,P)$ satisfies \eqref{NS} in the sense of distributions.\\
3.The function pair $(u,P)$ satisfies the following
local energy inequality:
\begin{align}\label{1.16}
&\int_{B_1(x_0)}|u(t,x)|^2\phi dx+2\int_{t_0-1}^t\int_{B_1(x_0)}|\nabla u|^2\phi dxds\\
\leq &\int_{t_0-1}^t\int_{B_1(x_0)}|u|^2(\partial_t\phi+\nu\Delta\phi)+(|u|^2+2P)u\cdot\nabla\phi dxds.\notag
\end{align} for every nonnegative $\phi\in C^\infty_0(Q_1(z_0))$.
\end{definition}
We say a point $z_0$ is a regular point of a solution $u$ to \eqref{NS} if there exists a non-empty neighborhood $\mathcal{O}_{z_0}$ of $z_0$ such that $u\in L^\infty(\mathcal{O}_{z_0})$. The complement of the set of regular points will be called the singular set.

Assume $\Omega\subset\mathbb{R}^3$. We use $L^q((0,T]; L^p(\Omega))$ to denote the space of measurable functions with the following norm
\begin{align*}
||f||_{L^q([0,T];L^p(\Omega))}=\left\{\begin{array}{ll}
(\int_0^T(\int_{\Omega}|f(t,x)|^pdx)^{\frac{q}{p}}dt)^{\frac{1}{q}},~1\leq q<\infty,\\
ess\sup_{t\in(0,T]}||f(t,\cdot)||_{L^p(\Omega)},~q=\infty.
\end{array}\right.
\end{align*}

The Lorentz space $L^{r,s}([0,T])$ is the space of measurable functions with the following norm:
\begin{align*}
||f||_{L^{r,s}([0,T])}=\left\{\begin{array}{ll}
(\int_0^\infty\sigma^{s-1}|\{x\in[0,T]:|f(x)|>\sigma\}|^{\frac{s}{r}}d\sigma)^\frac{1}{s},~1\leq s<\infty,\\
\sup_{\sigma>0}\sigma|\{x\in[0,T]:|f(x)|>\sigma\}|^\frac{1}{r},~s=\infty.
\end{array}\right.
\end{align*}

Let $(u, P)$ be a solution of \eqref{NS}. Introduce the scaling
\begin{align}\label{1.17}
u_\lambda(t,x)=\lambda u(\lambda^2t,\lambda x);~~P_\lambda(t,x)=\lambda^2P(\lambda^2t,\lambda x),
\end{align}for arbitrary $\lambda>0$. Then the function pair $(u_\lambda,P_\lambda)$ is also a solution of \eqref{NS}.

 We introduce
the following invariant quantities, which are invariant under the natural scaling \eqref{1.17}:
\begin{align*}
A(u,r,z)=\sup_{t-r^2\leq s\leq t}\frac{1}{r}\int_{B_r(x)\times \{s\}}|u|^2dx;~ B(u,r,z)=\frac{1}{r}\int\int_{Q_r(z)}|\nabla u|^2dxdt,\\
C(u,r,z)=\frac{1}{r^2}\int\int_{Q_r(z)}|v|^3dxdt;~ D(P,r,z)=\frac{1}{r^2}\int\int_{Q_r(z)}|P|^{\frac{3}{2}}dxdt.
\end{align*}
For simplicity, we introduce the notations
\begin{align*}
A(u,r)=A(u,r,0);~B(u,r)=&B(u,r,0);~C(u,r)=C(u,r,0);~D(P,r)=D(P,r,0).
\end{align*}

Throughout this paper, $u_{x_0,\rho}\doteq\frac{1}{|B_\rho|}\int_{B_\rho(x_0)}udx$ and $C$ denotes an absolute and often large positive number which can change from line to line.

\section{The proofs of Main results}
We first show some crucial lemmas.
\begin{lemma}\label{l1}
Let $z_0=(x_0,t_0)$ and $\frac{1}{q}+\frac{1}{p}=\frac{1}{2},~2\leq p\leq\infty.$  Assume $(u,P)$ be a suitable weak solution to \eqref{NS} on $Q_1(z_0)$ satisfying
\begin{align}\label{ass}
||\sup_{\rho\leq1}(\frac{1}{\rho}\int_{B_\rho(x_0)}|u(t,x)-u(t)_{x_0,\rho}|^p dx)^\frac{1}{p}||_{L^{q,\infty}([t_0-1,t_0])}=&M<\infty\\
or~~~~~~~~||\sup_{\rho\leq1}(\frac{1}{\rho}\int_{B_\rho(x_0)}|u(t,x)|^p dx)^\frac{1}{p}||_{L^{q,\infty}([t_0-1,t_0])}=&M<\infty.\end{align}Then,\\ if $2\leq p<3$, it holds
\begin{align}\label{2.2}
C(u,r,z_0)\leq C\frac{r}{\rho}C(u,\rho,z_0)+C(\frac{\rho}{r})^2B(u,\rho,z_0)^{\frac{9-3p}{6-p}}M^\frac{3p}{6-p}.
\end{align}
if $3\leq p\leq 6$, it holds
\begin{align}\label{2.3}
C(u,r,z_0)\leq C\frac{r}{\rho}C(u,\rho,z_0)+C(\frac{\rho}{r})A(u,\rho,z_0)^\frac{p-3}{p-2}M^\frac{p}{p-2},
\end{align}
if $6<p\leq\infty$, it holds
\begin{align}\label{2.4}
C(u,r,z_0)\leq C\frac{r}{\rho}C(u,\rho,z_0)+C(\frac{\rho}{r})^\frac{3}{2}A(u,\rho,z_0)^\frac{3}{4}M^\frac{3}{2}
\end{align}where $C$ is a positive absolute constant independent on $p$.
\end{lemma}

\textit{\bf Proof}\quad We first consider that the assumption (2.1) is holding. Let $r<\rho\leq 1$ and define $f_p(t)=(\sup_{\rho\leq1}\frac{1}{\rho}\int_{B_\rho(x_0)}|u(t,x)-u_{x_0,\rho}|^pdx)^\frac{1}{p}$.  At almost every time $t\in(t_0-\rho^2,t_0]$ we estimate
\begin{align}\label{2.5}
\int_{B_r(x_0)}|u|^3dx\leq C|B_r||u_{x_0,\rho}|^3+C\int_{B_r(x)}|u-u_{x_0,\rho}|^3dx
=I_1+I_2.
\end{align}

For $I_1$, we have
\begin{align}\label{2.6}
I_1=C|B_r(x_0)||\frac{1}{|B_\rho(x_0)|}\int_{B_\rho(x_0)}udy|^3
\leq C(\frac{r}{\rho})^3\int_{B_\rho(x_0)}|u|^3dx.
\end{align}

We now estimate $I_2$.

If $2< p<3$, we estimate $I_2$ as follows
\begin{align*}
I_2\leq C||u-u_{x_0,\rho}||^{\frac{3p}{6-p}}_{L^p}||u-u_{x_0,\rho}||^{\frac{6(3-p)}{6-p}}_{L^6}.
\end{align*} Integrating with respect to time from $t_0-r^2$ to $t_0$ and using H\"{o}lder's inequality, we obtain
\begin{align}\label{2.7}
&\int_{t_0-r^2}^{t_0}\int_{B_r(x_0)}|u-u_{x_0,\rho}|^3dxds\\\leq &(\int_{t_0-r^2}^{t_0}||u-u_{x_0,\rho}||^2_{L^6(B_\rho(x_0))}ds)^{\frac{9-3p}{6-p}}(\int_{t_0-r^2}^{t_0}||u-u_{x_0,\rho}||^{\frac{3p}{2p-3}}_{L^p(B_\rho(x_0))}ds)^{\frac{2p-3}{6-p}}\notag\\
\leq &(\int_{t_0-r^2}^{t_0}\int_{B_\rho(x_0)}|\nabla u|^2dxds)^{\frac{9-3p}{6-p}}\rho^{\frac{3}{6-p}}(\int_{t_0-r^2}^{t_0}f^\frac{3p}{2p-3}_p(s)ds)^{\frac{2p-3}{6-p}}\notag.
\end{align}
By using the assumption $||f_p||_{L^{q,\infty([t_0-1,t_0])}}=M$ and $\frac{1}{q}+\frac{1}{p}=\frac{1}{2}$, we have
\begin{align}\label{2.8}
&\int_{t_0-r^2}^{t_0}f_p(s)^{\frac{3p}{2p-3}}ds\\
=&\frac{3p}{2p-3}\int_{0}^{\infty}\sigma^\frac{3+p}{2p-3}|\{s\in[t-r^2,t];f_p(s)>\sigma\}|d\sigma\notag\\
\leq&\frac{3p}{2p-3}\{\int_{0}^{R}\sigma^\frac{3+p}{2p-3}r^2d\sigma+M^{\frac{2p}{p-2}}\int_{R}^{\infty}\sigma^{\frac{3+p}{2p-3}-\frac{2p}{p-2}}d\sigma\}\notag\\
\leq &R^\frac{3p}{2p-3}r^2+(3-\frac{6}{p})R^{\frac{3p}{2p-3}-\frac{2p}{p-2}}M^\frac{2p}{p-2}\notag\\
\leq &(4-\frac{6}{p})r^{\frac{p}{2p-3}}M^{\frac{3p}{2p-3}},\notag
\end{align}where  we take $R=r^{-\frac{p-2}{p}}M$.

When $p=2$, in the estimate \eqref{2.8}, we choose $R=M$ and obtain
\begin{align*}
\int_{t_0-r^2}^{t_0}f_p(s)^{\frac{3p}{2p-3}}ds=6\int_0^M\sigma^5r^2d\sigma=r^2M^6.
\end{align*} This means that the conclusion in \eqref{2.8} is still holding for $p=2$.

Substituting \eqref{2.8} into \eqref{2.7} implies
\begin{align}\label{2.9}
&\int_{t_0-r^2}^{t_0}\int_{B_r(x_0)}|u-u_{x_0,\rho}|^3dxds\\\leq &(4-\frac{6}{p})^{\frac{2p-3}{6-p}}\rho^{\frac{3}{6-p}}r^{\frac{p}{6-p}}(\int_{t_0-\rho^2}^{t_0}\int_{B_\rho(x_0)}|\nabla u|^2dxds)^{\frac{9-3p}{6-p}}M^\frac{3p}{6-p}.\notag
\end{align} Combining \eqref{2.9} with \eqref{2.5}-\eqref{2.6}, we get
\begin{align*}
&\int_{Q_r(z_0)}|u|^3dxdt\\\leq& C(\frac{r}{\rho})^3\int_{Q_\rho(z_0)}|u|^3dxds+C\rho^\frac{3}{6-p}r^\frac{p}{6-p}(\int_{Q_{\rho}(z_0)}|\nabla u|^2dxds)^\frac{9-3p}{6-p}M^\frac{3p}{6-p}\\
\leq &(\frac{r}{\rho})^3\int_{Q_\rho(z_0)}|u|^3dxds+C\rho^{2-\frac{p}{6-p}}r^\frac{p}{6-p}B(u,\rho,z_0)^\frac{9-3p}{6-p}M^\frac{3p}{6-p}
\end{align*}where we have used the fact $(4-\frac{6}{p})^\frac{2p-3}{6-p}\leq 4$ for $2\leq p<3$.
Multiplying this estimate by $\frac{1}{r^2}$, we obtain \eqref{2.2}.

 If $3\leq p\leq6$, we deduce, using interpolation inequality
\begin{align}\label{2.10}
I_2\leq& C(\int_{B_\rho(x_0)}|u-u_{x_0,\rho}|^2dx)^\frac{p-3}{p-2}(\int_{B_\rho(x_0)}|u-u_{x_0,\rho}|^pdx)^\frac{1}{p-2}\\
\leq & C\rho A(u,\rho,z_0)^{\frac{p-3}{p-2}}(\frac{1}{\rho}\int_{B_\rho(x_0)}|u-u_{x_0,\rho}|^pdx)^\frac{1}{p-2}.\notag
\end{align}Summing up the estimates for $I_1$ and $I_2$ and integrating with respect to time from $t_0-r^2$ to $t_0$, we obtain
\begin{align}\label{2.11}
&\int_{t_0-r^2}^{t_0}\int_{B_r(x_0)}|u|^3dxds\\\leq &C(\frac{r}{\rho})^3\int_{t_0-\rho^2}^{t_0}\int_{B_\rho(x_0)}|u|^3dxds+C \rho A(u,\rho,z_0)^\frac{p-3}{p-2}\int_{t_0-r^2}^{t_0}f^\frac{p}{p-2}_p(s)ds.\notag
\end{align}By the assumptions, we obtain  $||f_p||_{L^{q,\infty}[t_0-1,t_0]}=M$ with $\frac{1}{q}+\frac{1}{p}=\frac{1}{2}$. It follows
\begin{align}\label{2.12}
&\int_{t_0-r^2}^{t_0}f^\frac{p}{p-2}_p(s)ds\\
=&\frac{p}{p-2}\int_{0}^{\infty}\sigma^\frac{2}{p-2}|\{s\in[t_0-r^2,t_0]: f_p(s)
>\sigma\}|d\sigma\notag\\
=&\frac{p}{p-2}\{\int_{0}^{R}\sigma^\frac{2}{p-2}|\{s\in[t_0-r^2,t_0]: f_p(s)
>\sigma\}|d\sigma\notag\\&+
\int_{R}^{\infty}\sigma^\frac{2}{p-2}|\{s\in[t_0-r^2,t_0]: f_p(s)
>\sigma\}|d\sigma\notag\\
\leq &\frac{p}{p-2}\int_{0}^{R}\sigma^\frac{2}{p-2}r^2d\sigma+\frac{p}{p-2}\int_{R}^{\infty}\sigma^{\frac{2}{p-2}-\frac{2p}{p-2}}d\sigma M^{\frac{2p}{p-2}}\notag\\
=&r^2R^\frac{p}{p-2}+R^\frac{-p}{p-2}M^\frac{2p}{p-2}\notag\\
=&2rM^\frac{p}{p-2}\notag
\end{align}where we choose $R=r^{-\frac{p-2}{p}}M$. Substituting \eqref{2.12} into \eqref{2.11}, it follows
\begin{align}\label{2.13}
&\int_{t_0-r^2}^{t_0}\int_{B_r(x_0)}|u|^3dxdt\\\leq &C(\frac{r}{\rho})^3\int_{t_0-\rho^2}^{t_0}\int_{B_\rho(x_0)}|u|^3dxdt+C r\rho A(u,\rho,z_0)^\frac{p-3}{p-2}M^\frac{p}{p-2}.\notag
\end{align} Multiplying \eqref{2.13} by $\frac{1}{r^2}$, we get
\begin{align*}
\frac{1}{r^2}\int_{Q_r(z)}|u|^3 dxdt\leq C(\frac{r}{\rho})\frac{1}{\rho^2}\int_{Q_\rho(z)}|u|^3dxdt+C(\frac{\rho}{r})A(u,\rho,z_0)^\frac{p-3}{p-2}M^\frac{p}{p-2}.
\end{align*}
This means \eqref{2.3}.

If $6<p\leq\infty$, we estimate $I_2$ by using H\"{o}lder's inequality as follows
\begin{align}\label{2.14}
I_2\leq&\int_{B_\rho(x_0)}|u-u_{x_0,\rho}|^\frac{3}{2}|u-u_{x_0,\rho}|^\frac{3}{2}dx\\
\leq&C(\int_{B_\rho(x_0)}|u-u_{x_0,\rho}|^2dx)^\frac{3}{4}(\int_{B_\rho(x_0)}|u-u_{x_0,\rho}|^p)^\frac{3}{2p}\rho^\frac{3(p-6)}{4p}\notag\\
\leq& C\rho^{\frac{3}{4}+\frac{3}{2p}+\frac{3(p-6)}{4p}}A(u,\rho,z_0)^\frac{3}{4}f^\frac{3}{2}_p(s)ds.\notag
\end{align}
Summing up the estimates for $I_1$ and $I_2$ and integrating with respect to time from $t_0-r^2$ to $t_0$, we obtain
\begin{align}\label{2.15}
&\int_{t_0-r^2}^{t_0}\int_{B_r(x_0)}|u|^3dxds\\\leq &C(\frac{r}{\rho})^3\int_{t_0-\rho^2}^{t_0}\int_{B_\rho(x_0)}|u|^3dxds+C \rho^\frac{3p-6}{2p} A(u,\rho,z_0)^\frac{3}{4}\int_{t_0-r^2}^{t_0}f^\frac{3}{2}_p(s)ds.\notag
\end{align} Using the similar estimates for \eqref{2.8} or \eqref{2.12}, we obtain
\begin{align}\label{2.16}
\int_{t_0-r^2}^{t_0}f^\frac{3}{2}_p(s)ds\leq(1+3\frac{p-2}{p+6}) r^{2-\frac{3p-6}{2p}}M^\frac{3}{2}.
\end{align}Substituting \eqref{2.16} into \eqref{2.15} implies
\begin{align}\label{2.17}
&\frac{1}{r^2}\int_{t_0-r^2}^{t_0}\int_{B_r(x_0)}|u|^3dxds\\\leq& C\frac{r}{\rho}\frac{1}{\rho^2}\int_{t_0-\rho^2}^{t_0}\int_{B_\rho(x_0)}|u|^3dxds+C(\frac{\rho}{r})^\frac{3p-6}{2p}A(u,\rho,z_0)^\frac{3}{4}M^\frac{3}{2}\notag\\
\leq &C\frac{r}{\rho}C(u,\rho,z_0)+C(\frac{\rho}{r})^\frac{3}{2}A(u,\rho,z_0)^\frac{3}{4}M^\frac{3}{2}\notag
\end{align}where we have used the facts $1+3\frac{p-2}{p+6}\leq4$ and $(\frac{\rho}{r})^\frac{3p-6}{2p}\leq(\frac{\rho}{r})^\frac{3}{2}$. We thus show \eqref{2.4} and complete the proof of Lemma
\ref{l1} under the assumption (2.1).

If the assumption $(2.2)$ is holding, we denote $f_p(t)=\sup_{\rho\leq1}(\frac{1}{\rho}\int_{B_\rho(x_0)}|u(t,x)|^pdx)^\frac{1}{p}$ and modify the processes of proofs as follows.

In the case $2\leq p<3$, we can replace $||u-u_\rho||_{L^p(B_\rho)}$ by $C_p||u||_{L^p(B_\rho)}$ in \eqref{2.7} and repeat the processes of proofs for \eqref{2.8}-\eqref{2.9} to get \eqref{2.2}. The difference is that in this case, the constant $C_p$ is depended on $p$. Noticing $2\leq p<3$, we can choose a large enough constant $C$ to get rid of the dependence on $p$.

If $3\leq p\leq 6$ or $6<p\leq\infty$, we just need to replace $I_2=\int_{B_r(x_0)}|u-u_{x_0,\rho}|^3dx$ by  $C\int_{B_\rho(x_0)}|u|^3dx$ and repeat the processes of proofs step by step to get \eqref{2.3} and \eqref{2.4}.

\begin{lemma}\label{l2}
Let $z_0=(x_0,t_0)$ and $\frac{1}{q}+\frac{1}{p}=\frac{1}{2},~2\leq p\leq\infty.$  Assume $(u,P)$ be a suitable weak solution to \eqref{NS} on $Q_1(z_0)$ satisfying
\begin{align}\label{as1}
||\sup_{\rho\leq1}(\frac{1}{\rho}\int_{B_\rho(x_0)}|u(t,x)-u(t)_{x_0,\rho}|^p dx)^\frac{1}{p}||_{L^{q,\infty}([t_0-1,t_0])}=M<\infty,\\
or~~~~~~~||\sup_{\rho\leq1}(\frac{1}{\rho}\int_{B_\rho(x_0)}|u(t,x)|^p dx)^\frac{1}{p}||_{L^{q,\infty}([t_0-1,t_0])}=M<\infty,\end{align}then there exists a constant $\rho_0>0$ only depended on $A(u,1,z_0)$, $B(u,1,z_0)$, $C(u,1,z_0)$ and $D(P,1,z_0)$, such that for $r\leq\rho_0$, it follows
\begin{align}\label{2.19}
A(u,r,z)+B(u,r,z)+C(u,r,z)+D(P,r,z)\leq C(M).
\end{align}
\end{lemma}

\textit{\bf Proof}\quad
Without loss of generality, we set $z_0=0$.
Let $\phi(t,x)=\chi(t,x)\psi(t,x)$ where $\chi$ is cut-off function which equals 1 in
$Q_{\frac{1}{2}\rho}$ and vanishes outside of $Q_{\frac{3}{4}\rho}$. Then let $\psi=(4\pi(r^2-t))^{-\frac{3}{2}}e^{-\frac{|x|^2}{4(r^2-t)}}$. Direct computations show that $\phi\geq0$ and
\begin{align*}
\partial_t\phi+\triangle\phi=&0~in~Q_{\frac{1}{2}\rho},\\
|\partial_t\phi+\triangle\phi|\leq& C\rho^{-5}~in~Q_{\rho},\\
C^{-1}r^3\leq \phi\leq Cr^{-3};~|\nabla\phi|\leq& Cr^{-4}~in~Q_r,\\
\phi\leq C\rho^{-3};~|\nabla\phi|\leq& C\rho^{-4}~in~Q_\rho-Q_{\frac{3}{4}\rho}.
\end{align*}
Using $\phi$ as a test function in the local energy inequality \eqref{1.16}, we obtain
\begin{align}\label{2.20}
A(u,r)+B(u,r)\leq& C(\frac{r}{\rho})^2A(u,\rho)+C(\frac{\rho}{r})^2C(u,\rho)+C(\frac{\rho}{r})^2C^\frac{1}{3}(u,\rho)D^\frac{2}{3}(P,\rho)\\
\leq& C(\frac{r}{\rho})^2A(u,\rho)+C(\frac{\rho}{r})^2C(u,\rho)+C(\frac{\rho}{r})^2D(P,\rho)\notag.
\end{align}

We now show some bounds on $D(u,r)$. Let $\eta(x)$ be a cut-off function which equals 1 in $B_{\frac{3\rho}{4}}$ and vanishes outside of $B_\rho$. Let $P_1$ satisfy $-\Delta P_1=\partial_{x_i}\partial_{x_j}(u_iu_j\eta)$ and $P_2=P-P_1$. Then, it is clear that $\Delta P_2=0$ in $B_\frac{3\rho}{4}$. By using the Calder\'{o}n-Zygmund inequality, we have
\begin{align*}
\int_{B_\rho}|P_1|^{\frac{3}{2}}dx\leq C(\int_{B_\rho}|u|^3dx).
\end{align*}By the properties of the harmonic functions, we infer that for $r\leq\frac{\rho}{2}$,
\begin{align*}
\int_{B_r}|P_2|^{\frac{3}{2}}dx\leq Cr^3\sup_{x\in B_r}|P_2(x)|^\frac{3}{2}\leq C(\frac{r}{\rho})^3\int_{B_\rho}|P_2|^{\frac{3}{2}}dx.
\end{align*}It then follows that for $0<r\leq\frac{\rho}{2}$
\begin{align*}
&\int_{B_r}|P|^\frac{3}{2}dx\\\leq& C(\int_{B_\rho}|u|^3dx)+C(\frac{r}{\rho})^3\int_{B_\rho}|P-P_1|^\frac{3}{2}dx\\
\leq& C(\int_{B_\rho}|u|^3dx)+C(\frac{r}{\rho})^3\int_{B_\rho}|P|^\frac{3}{2}dx.
\end{align*}Integrating with respect to t from $-r^2$ to 0, we obtain, using H\"{o}lder inequality,
\begin{align*}
\int_{Q_r}|P|^{\frac{3}{2}}dxdt\leq C\int_{Q_\rho}|u^3|dxdt+C(\frac{r}{\rho})^3\int_{Q_\rho}|P|^\frac{3}{2}dxdt.
\end{align*} This implies
\begin{align}\label{2.21}
D(P,r)\leq C\frac{r}{\rho}D(P,\rho)+C(\frac{\rho}{r})^2C(u,\rho).
\end{align}

We now show some crucial bounds for $C(u,r)$.

{\bf In the case $2\leq p<3$}\quad Noticing \eqref{2.2},
 we have by using Young's inequality
\begin{align}\label{2.22}
C^\frac{7}{6}(u,r)\leq& C(\frac{r}{\rho})^\frac{6}{7}C^\frac{7}{6}(u,\rho)+C(\frac{\rho}{r})^\frac{7}{3}B^\frac{7(3-p)}{2(6-p)}(u,\rho)M^\frac{7p}{2(6-p)}\\\leq &C(\frac{r}{\rho})^\frac{6}{7}C^\frac{7}{6}(u,\rho)+C(\frac{\rho}{r})^{\frac{7}{3}+\frac{7}{6}\frac{7(3-p)}{2(6-p)}}((\frac{r}{\rho})^\frac{7}{6}B(u,\rho))^\frac{7(3-p)}{2(6-p)}M^\frac{7p}{2(6-p)}\notag\\
\leq &C(\frac{r}{\rho})^\frac{6}{7}C^\frac{7}{6}(u,\rho)+C(\frac{\rho}{r})^5((\frac{r}{\rho})^\frac{7}{6}B(u,\rho))^\frac{7(3-p)}{2(6-p)}M^\frac{7p}{2(6-p)},\notag\\
\leq &C(\frac{r}{\rho})^\frac{6}{7}C^\frac{7}{6}(u,\rho)+C(\frac{r}{\rho})^\frac{7}{6}B(u,\rho)+C(\frac{\rho}{r})^{40}M^\frac{7p}{5p-9},\notag
\end{align}where we have used the facts $\frac{\rho}{r}>1$ and $2\leq p<3$.

By using Young's inequality, we deduce from \eqref{2.20}
\begin{align}\label{2.23}
&A(u,r)+B(u,r)\\\leq &C(\frac{r}{\rho})^2A(u,\rho)+(\frac{r}{\rho})^\frac{7}{6}C^\frac{7}{6}(u,\rho)
+(\frac{r}{\rho})^\frac{8}{7}D^\frac{8}{7}(P,\rho)+C((\frac{\rho}{r})^{21}+(\frac{\rho}{r})^{24}).\notag
\end{align}
Similarly, we obtain by using \eqref{2.21}
\begin{align}\label{2.24}
D(P,r)^\frac{8}{7}\leq& C(\frac{r}{\rho})^\frac{8}{7}D^\frac{8}{7}{(P,\rho)}+C(\frac{\rho}{r})^\frac{16}{7}C(u,\rho)^\frac{8}{7}\\
\leq& C(\frac{r}{\rho})^\frac{8}{7}D^\frac{8}{7}{(P,\rho)}+(\frac{r}{\rho})^\frac{7}{6}C(u,\rho)^\frac{7}{6}
+C(\frac{\rho}{r})^{168}\notag.\end{align}

Define $G(r)\equiv A(u,r)+B(u,r)+C^\frac{7}{6}(u,r)+D^\frac{8}{7}(P,r)$. Summing up the estimates \eqref{2.22}-\eqref{2.24} implies
\begin{align}\label{2.25}
G(r)\leq C(\frac{r}{\rho})^\frac{8}{7}G(\rho)+C(1+M^{\frac{7p}{5p-9}})(\frac{\rho}{r})^{168}
\end{align} where we have used the fact $\frac{r}{\rho}<1$.

 Fix $\theta=\min\{\frac{1}{2},\frac{1}{C^7}\}$ and set $r=\theta^{k}\rho$ for $k\in\mathbb{N}$. \eqref{2.25} yields
\begin{align}\label{2.26}
G(\theta^{k}\rho)\leq\theta G(\theta^{k-1}\rho)+C(1+M^{\frac{7p}{5p-9}})\theta^{-168}.
\end{align}By a standard iterative argument, we deduce that
\begin{align}\label{2.27}
G(r)\leq \frac{r}{\rho}G(\rho)+C(1+M^{\frac{7p}{5p-9}})~for~r\leq\rho\leq 1.
\end{align}We now first take $\rho=1$ then choose $\rho_{01}$ satisfying $\frac{\rho_{01}}{1}G(1)\leq 1$, it follows
\begin{align}\label{2.28}
G(r)\leq C(M^\frac{7p}{5p-9})~for~r\leq\rho_{01}.
\end{align}

{\bf In the case $3\leq p\leq 6$}\quad From \eqref{2.3}, it is clear that
\begin{align}\label{2.29}
&C(u,r)^\frac{7}{6}\\\leq& C(\frac{r}{\rho})^\frac{7}{6}C(u,\rho)^\frac{7}{6}+C (
\frac{\rho}{r})^\frac{7}{6}A(u,\rho)^{\frac{7}{6}\frac{p-3}{p-2}}M^{\frac{7}{6}\frac{p}{p-2}}\notag\\
\leq &C(\frac{r}{\rho})^\frac{7}{6}C(u,\rho)^\frac{7}{6}+C (
\frac{\rho}{r})^{\frac{7}{6}+(\frac{7}{6})^2\frac{p-3}{p-2}}((\frac{r}{\rho})^\frac{7}{6}A(u,\rho))^{\frac{7}{6}\frac{p-3}{p-2}}M^{\frac{7}{6}\frac{p}{p-2}}\notag\\
\leq & C(\frac{r}{\rho})^\frac{7}{6}C(u,\rho)^\frac{7}{6}+(\frac{r}{\rho})^\frac{7}{6}A(u,\rho)+C (\frac{\rho}{r})^{\frac{6(p-2)}{9-p}(\frac{7}{6}+\frac{49}{36}\frac{p-3}{p-2})}M^{\frac{7p}{9-p}}\notag\\
\leq &C(\frac{r}{\rho})^\frac{7}{6}C(u,\rho)^\frac{7}{6}+(\frac{r}{\rho})^\frac{7}{6}A(u,\rho)+C(\frac{\rho}{r})^{24}M^{\frac{7p}{9-p}}\notag
\end{align}where we have used the fact ${\frac{6(p-2)}{9-p}(\frac{7}{6}+\frac{49}{36}\frac{p-3}{p-2})}\leq 24$ for $3\leq p\leq 6$.
Collecting \eqref{2.23}-\eqref{2.24} and \eqref{2.29} implies
\begin{align}\label{2.30}
G(r)\leq C(\frac{r}{\rho})^\frac{8}{7}G(\rho)+C(1+M^\frac{7p}{9-p})(\frac{\rho}{r})^{168}.
\end{align}By using the similar computations in the estimates for \eqref{2.25}-\eqref{2.28}, we show that there exists a constant $\rho_{02}$ such that for $r\leq\rho_{02}$ it follows
\begin{align}\label{2.31}
G(r)\leq C(M^\frac{7p}{9-p}).
\end{align}

{\bf In the case $6<p\leq\infty$:}\quad
From \eqref{2.4} and Young's inequality, it is clear that
\begin{align}\label{2.32}
C(u,r)^\frac{7}{6}\leq C(\frac{r}{\rho})^\frac{7}{6}C(u,\rho)^\frac{7}{6}+(\frac{r}{\rho})^\frac{7}{6}A(u,\rho)+C(\frac{\rho}{r})^{23}M^{14}.
\end{align} Collecting \eqref{2.23}-\eqref{2.24} and \eqref{2.32} yieds
\begin{align}\label{2.33}
G(r)\leq C(\frac{r}{\rho})^\frac{8}{7}G(\rho)+C(\frac{\rho}{r})^{168}(M^{14}+1).
\end{align} By using similar computations in the estimates for \eqref{2.25}-\eqref{2.28}, we get that there exists a constant $\rho_{03}$ such that for $r\leq\rho_{03}$, it follows
\begin{align}\label{2.34}
G(r)\leq C(M^{14}).
\end{align}Collecting \eqref{2.28}, \eqref{2.31} and \eqref{2.34} and taking $\rho_0=\min\{\rho_{01},\rho_{02},\rho_{03}\}$, we thus obtain \eqref{2.19}.

\begin{lemma}\label{l3}
Let $z_0=(x_0,t_0)$ and $\frac{1}{q}+\frac{1}{p}=\frac{1}{2}$ with $2\leq p\leq\infty$. Assume $(u,P)$ be a suitable weak solution of \eqref{NS} in $Q_1(z_0)$. For any fixed $\varepsilon>0$, there exists two constants $\delta$ and $r^*$ depended on $\varepsilon$ such that if
\begin{align}\label{2.35}
||(\sup_{\rho\leq1}\frac{1}{\rho}\int_{B_\rho(x_0)}|u(t,x)-u(t)_{x_0,\rho}|^p dx)^\frac{1}{p}||_{L^{q,\infty}([t_0-1,t_0])}\leq \delta\\
or~~~~~~~||(\sup_{\rho\leq1}\frac{1}{\rho}\int_{B_\rho(x_0)}|u(t,x)|^p dx)^\frac{1}{p}||_{L^{q,\infty}([t_0-1,t_0])}\leq \delta\end{align}then it is holding
\begin{align}\label{2.36}
C(u,r^*,z_0)\leq\varepsilon,
\end{align}
\end{lemma}
\textit{\bf Proof}\quad Without loss of generality, we assume $z_0=0$ and $\delta\leq 1$. In view of Lemma \ref{l2}, we have that for $\rho\leq\rho_0$, it is holding
\begin{align}\label{2.37}
&C(u,\rho,z)+D(P,\rho,z)+A(u,\rho,z)+B(u,\rho,z)\leq C
\end{align}where $C$ is an absolute constant.

 If $2\leq p<3$,  we deduce by choosing $\rho= \rho_0$ in \eqref{2.2} and using \eqref{2.37}
\begin{align}\label{2.38}
C(u,r)\leq& \frac{r}{\rho_0}C+C(\frac{\rho_0}{r})^2C^\frac{9-3p}{6-p}\delta^\frac{3p}{6-p}\\
\leq &\frac{r}{\rho_0}C+C^2(\frac{\rho_0}{r})^2\delta^\frac{3}{2}\notag
\end{align}where we have used $\frac{9-3p}{6-p}\leq 1$ and $\frac{3p}{6-p}\geq\frac{3}{2}$. In \eqref{2.38}, we first choose $r^*=\frac{\varepsilon\rho_0}{2C}$ then take
 $\delta_1\leq \frac{\varepsilon^2}{4C^\frac{8}{3}}$, it follows that
 \begin{align}\label{2.39}
 C(u,r^*)\leq\varepsilon.
 \end{align}

 If $3\leq p\leq6$, by choosing $\rho=\rho_0$ in \eqref{2.3} and using \eqref{2.37}, we obtain also that
 \begin{align}\label{2.40}
 C(u,r)\leq& C\frac{r}{\rho_0}C(u,\rho_0)+C(\frac{\rho_0}{r})A(u,\rho)^\frac{p-3}{p-2}\delta^\frac{p}{p-2}\\
 \leq &\frac{r}{\rho_0}C+C(\frac{\rho_0}{r}) C\delta^\frac{3}{2}\notag
 \end{align}where we have used $\frac{p-3}{p-2}\leq 1$ and $\frac{p}{p-2}\geq\frac{3}{2}$ for $3\leq p\leq 6$.
 We now first choose $r^*=\frac{\varepsilon\rho_0}{2C}$ then take $\delta_2\leq\frac{\varepsilon^\frac{4}{3}}{4^\frac{2}{3}C^2}$, it follows \eqref{2.39} again.

If $6<p\leq\infty$, by choosing $\rho=\rho_0$ in \eqref{2.4} and using \eqref{2.37}, we obtain by similar computations that
\begin{align}\label{2.41}
C(u,r)\leq C\frac{r}{\rho_0}+C^2
(\frac{\rho_0}{r})^\frac{3}{2}\delta^\frac{3}{2}.\end{align}
We now first choose $r^*=\frac{\varepsilon\rho_0}{2C}$ then take $\delta_3\leq\frac{\varepsilon^\frac{5}{3}}{2^\frac{5}{3}C^\frac{7}{3}}$, it follows \eqref{2.39}.  Choosing $\delta=\min\{\delta_1,\delta_2,\delta_3\}$ yields Lemma \ref{l3}.

To get the concentration rate including both time scale and space scale, we need an embedding theorem from the Lorentz space $L^{p,\infty}$ to a Morrey type space
\begin{lemma}\label{l4}
For any given $r>0$ and $2\leq p\leq\infty$, it follows
\begin{align}\label{2.42}
(\sup_{\eta\leq r}\frac{1}{\eta}\int_{B_\eta}|u|^pdx)^\frac{1}{p}\leq C||u||_{L^{\frac{3p}{2},\infty}(B_r)}
\end{align}where $C>0$ is a constant independent on $p$.
\end{lemma}
\textit{\bf Proof}\quad This conclusion is a direct computation as follows
\begin{align}\label{2.43}
\int_{B_\eta}|u|^pdx=&p\int_{0}^{\infty}\sigma^{p-1}|\{x\in B_\eta:|u(x)|>\sigma\}|d\sigma\\
\leq&Cp[\int_{0}^{R}\sigma^{p-1}\eta^3d\sigma+\int_{R}^{\infty}\sigma^{p-1-\frac{3p}{2}}d\sigma||u||^\frac{3p}{2}_{L^{\frac{3p}{2},\infty}(B_\eta)}]\notag\\
\leq& C[R^p\eta^3+2R^\frac{-p}{2}||u||^\frac{3p}{2}_{L^{\frac{3p}{2},\infty}(B_\eta)}]\notag\\
\leq& C\eta||u||^p_{L^{\frac{3p}{2},\infty}(B_\eta)}\notag
\end{align}where we take $R=\eta^\frac{-2}{p}||u||_{L^{\frac{3p}{2},\infty}(B_\eta)}$. This yields \eqref{2.42}.

To prove the local regularity for the suitable weak solution to \eqref{NS}, we need a criterion for partial regularity due to Wolf \cite{Wolf}.
\begin{lemma}\label{l5}\cite{Wolf}
For every $3 \leq s, q\leq\infty$ there exists a constant $\varepsilon^*=\varepsilon(s,p)>0$
with the following property: Let $u$ be a suitable weak solution to the Navier-Stokes
equations \eqref{NS} in $Q_r(z)$. If
\begin{align*}
||u||_{L^q(([t-r^2,t);L^p(B_r(x)))}\leq r^{\frac{2}{q}+\frac{3}{p}-1}\varepsilon^*,
\end{align*}Then $u$ is  H\"{o}lder continuous on $Q_{\frac {r}{2}}(z)$.
\end{lemma}

 We now start the proofs of Theorem \ref{mainr} and Theorem \ref{mains}.\\
\textit{\bf Proof of Theorem \ref{mainr}}\quad
 By using Lemma \ref{l3}, we obtain that for the given $\varepsilon^*=\varepsilon(3,3)$ in Lemma \ref{l5}, there exist two positive constants $r^*$ and $\delta$ such that if
 \begin{align*}
||\sup_{\eta\leq 1}(\frac{1}{\eta}\int_{B_\eta(x_0)}|u(x,t)-u_{x_0,\eta}|^pdx)^\frac{1}{p}||_{L^{q,\infty}[t_0-1,t_0]}\leq&\delta\\
or~~~~~~~ ||\sup_{\eta\leq 1}(\frac{1}{\eta}\int_{B_\eta(x_0)}|u(x,t)|^pdx)^\frac{1}{p}||_{L^{q,\infty}[t_0-1,t_0]}\leq&\delta
\end{align*} then it follows
\begin{align*}
C(u,r^*,z)\leq \varepsilon^*.
\end{align*} In view of Lemma \ref{l5} for $q=p=3$, we deduce that $z_0$ is a regular point. This yields Theorem \ref{mainr}.\\
\textit{\bf Proof of Theorem \ref{mains}}\quad Without loss of generality, we assume $0<r<1$. In the case $3<p\leq\infty$, if Theorem \ref{mains} is false, then there exists some $0<r_0<1$ such that for some $3< p_0\leq\infty$, $2\leq\nu_0\leq\frac{2p_0}{3}$ and $\frac{1}{\mu_0}+\frac{1}{\nu_0}=\frac{1}{2}$, it holds
\begin{align}\label{2.46}
\limsup_{t\to t_0}(t_0-t)^\frac{1}{\mu_0}r^{\frac{2}{\nu_0}-\frac{3}{p_0}}_0||u(t)||_{L^{p_0,\infty}(B_{r_0}(x_0))}\leq \delta^*.
\end{align} By using H\"{o}lder's inequality and Lemma \ref{l4}, we get
\begin{align}\label{2.48}
&\limsup_{t\to t_0}(t_0-t)^\frac{1}{\mu_0}\sup_{\eta<r_0}(\frac{1}{\eta}\int_{B_\eta(x_0)}|u|^{\nu_0} dx)^\frac{1}{\nu_0}\\
\leq& \limsup_{t\to t_0}(t_0-t)^\frac{1}{\mu_0}\sup_{\eta<r_0}(\eta^{2-\frac{9\nu_0}{2p_0}+\frac{3\nu_0}{2p}}(\frac{1}{\eta}\int_{B_\eta(x_0)}|u|^\frac{2p_0}{3} dx)^\frac{3\nu_0}{2p_0})^\frac{1}{\nu_0}\notag\\
\leq& \limsup_{t\to t_0}(t_0-t)^\frac{1}{\mu_0}\sup_{\eta<r_0}(\eta^{\frac{2}{\nu_0}-\frac{3}{p_0}}(\frac{1}{\eta}\int_{B_\eta(x_0)}|u|^\frac{2p_0}{3} dx)^\frac{3}{2p_0})\notag\\
\leq&  C\limsup_{t\to t_0}(t_0-t)^\frac{1}{\mu_0}\sup_{\eta<r_0}(\eta^{\frac{2}{\nu_0}-\frac{3}{p_0}}||u(t)||_{L^{p_0,\infty}(B_\eta(x_0)})\notag\\
\leq& C \limsup_{t\to t_0}(t_0-t)^\frac{1}{\mu_0}{r_0}^{\frac{2}{\nu_0}-\frac{3}{p}}||u(t)||_{L^{p_0,\infty}(B_{r_0}(x_0)}\notag
\leq C\delta^*\notag\\
= &\delta,\notag
\end{align} where $\delta$ is the same constant in Theorem \ref{mainr} and we choose $\delta^*=\frac{\delta}{C}$. This yields $||\sup_{\eta\leq r^*_0}(\frac{1}{\eta}\int_{B_\eta(x_0)}|u(t)|^{\nu_0} dx)^\frac{1}{\nu_0}||_{L^{\mu_0,\infty}[t_0-(r^*_0)^2,t_0]}\leq\delta$ for some $r^*_0\leq r_0$. From \eqref{2.48} and Theorem \ref{mainr}, we deduce $z_0$ is a regular point. This is a contradiction.

In the case $p=3$, if Theorem \ref{mains} is false, then there exists some $0<r_0<1$ such that it holds that
\begin{align}\label{2.49}
\limsup_{t\to t_0^-}||u(t,x)-u_{r,x_0}(t)||_{L^{3,\infty}(B_r(x_0))}\leq\delta^*.
\end{align} By using the fact $\min_{c\in\mathbb{R}}\int_{B_\eta(x_0)}|u-c|^2dx=\int_{B_\eta(x_0)}|u-u_{x_0,\eta}|^2dx$ and Lemma \ref{l4}, we deduce
\begin{align}\label{2.50}
&\limsup_{t\to t_0}\sup_{\eta<r_0}(\frac{1}{\eta}\int_{B_\eta(x_0)}|u-u_{x_0,\eta}|^2 dx)^\frac{1}{2}\\
\leq &C \limsup_{t\to t_0}\sup_{\eta<r_0}(\eta^{-1}(\int_{B_\eta(x_0)}|u-u_{x_0,r}|^2dx)^\frac{1}{2})\notag\\
=&C \limsup_{t\to t_0}||u(t)-u(t)_{x_0,r}||_{L^{3,\infty}(B_{r_0}(x_0))}\notag\\
\leq &C\delta^*=\delta.\notag
\end{align}This yields $||\sup_{\eta\leq r^*_0}(\frac{1}{\eta}\int_{B_\eta(x_0)}|u(t)-u_{x_0,\eta}|^{2} dx)^\frac{1}{2}||_{L^{\infty}[t_0-(r^*_0)^2,t_0]}\leq\delta$. From \eqref{2.50} and Theorem \ref{mainr}, we deduce $z_0$ is a regular point. This is also a contradiction. We thus complete the proof of Theorem \ref{mains}.
\section*{Acknowledgments}

The authors are supported by the Construct Program of the Key Discipline in Hunan Province and NSFC Grant No. 11871209.

\end{document}